\documentclass{amsart}

\usepackage{cite}

\usepackage[left=1.5in, right=1.5in, top=1.5in, bottom=1.5in]{geometry}
\usepackage{hyperref}
\usepackage{float}
\usepackage{tgcursor}
\usepackage[french,british]{babel}
\usepackage{amsmath}
\usepackage[dvipsnames]{xcolor}
\usepackage{hyperref} 
\hypersetup{backref=true,       
    pagebackref=true,               
    hyperindex=true,                
    colorlinks=true,                
    breaklinks=true,                
    urlcolor= black,                
    linkcolor= black,                
    bookmarks=true,                 
    bookmarksopen=false,
    filecolor=black,
    citecolor=black,
    linkbordercolor=green
}
\usepackage{filecontents}
\usepackage[utf8]{inputenc}
\usepackage{graphicx}
\usepackage{graphics}
\newcommand{\Grad}{\nabla\!\!\!\!\nabla}

\numberwithin{equation}{section}

\usepackage{bm,mathrsfs,dsfont,amssymb} 
\usepackage{color,soul,hyperref}

\usepackage{mathrsfs}
\usepackage[all]{xy}
\usepackage{tikz}
\usetikzlibrary{matrix}
\usepackage{multicol}
\usepackage{braket}
\usepackage{wrapfig}
\usepackage{listings}
\tikzset{node distance=2cm, auto}
\usepackage{array}
\newcolumntype{C}{>$c<$}
\usepackage{graphicx}
\usepackage{color}
\usepackage{lipsum}
\usepackage{float}
\usepackage{ stmaryrd }
\usepackage{ dsfont }
\usepackage{graphicx}
\usepackage{caption}
\usepackage{subcaption}
\usepackage{booktabs}
\usepackage{amssymb}
\usepackage{amsthm}
\usepackage{mathtools}
\newtheorem{prop}{Proposition}

\newtheorem{lemma}{Lemme}
\newtheorem{definition}{définition}
\newtheorem{thm}{Théorème}
\newtheorem*{prL}{Preuve du lemme}
\newtheorem*{prP}{Preuve de la proposition}
\newtheorem*{prT}{Preuve du Théorème}

\usepackage{mathrsfs}

\begin{document}
\selectlanguage{french}
\title[]{Morse homology: orientation of the moduli space of gradient flow lines, coherence and applications}

\author{Mathieu Giroux}
\address{department of Physics, McGill University, 3600 rue University, Montréal, QC Canada H3A 2T8}
\email{girouxma@hep.physics.mcgill.ca}

\date{\today}

\begin{abstract}
\textcolor{white}{a}\\

In this paper, we shall compute the chain complex and the corresponding homology of some Morse function $f$ over integer coefficients. The definition of the correct boundary operator requires a careful construction of moduli space of (pseudo)gradient flow lines orientations. We will then apply this construction in the computation of these homology groups on 4-manifolds.

-----

Dans ce papier, nous calculerons le complexe de chaînes défini sur les points critiques d'une fonction de Morse $f$ et l'homologie correspondante dans
le cas des coefficients entiers. La définition de l'opérateur de bord demande de définir correctement les orientations
des espaces de modules de flots induites par le (pseudo)gradient de $f$. On appliquera alors cette construction au calcul explicite de ces homologies sur les variétés de
dimension 4.

\end{abstract}
\maketitle
\smallskip
\noindent \textbf{Mots clés.} Théorie de Morse, homologie, orientation, cohérence.

\tableofcontents

\newpage
\section{Introduction}
Un bon nombre de questions venant de diverses branches des mathématiques mène au problème d'analyser la topologie d'un complexe simplicial. Cela dit, il n'existe que peu de techniques générales disponibles pour nous aider dans une telle tâche. Nous pouvons néanmoins apprécier que certaines théories très générales aient été développées dans cette optique, du moins dans le cas des variétés lisses. L'une des théories les plus puissantes et les plus utiles dans ce contexte est la théorie de Morse -- voir \cite{frauenfelder2020moduli,milnor2016morse} pour une couverture presque complète des termes utilisés tout au long de ce papier. En dimension finie ou en dimension infinie, cette théorie joue un rôle crucial dans la recherche
mathématique internationale actuelle \cite{grunert2019pl, laudenbach2020poincarlefschetz, Papageorgiou_2020, magao2020eternal}. En dimension finie (cas classique), cette homologie permet au moyen d'une fonction
de Morse réelle, c'est-à-dire d'une fonction réelle non dégénérée en ses points critiques (hessien non dégénéré)
définie sur une variété différentiable quelconque, de décomposer la variété en morceaux élémentaires qui
permettent de construire explicitement l'homologie de la variété. La version de Witten \cite{witten1982supersymmetry} permet de calculer cette
homologie en construisant un complexe de chaînes défini sur les points critiques de la fonction $f$ et en
définissant l'opérateur bord par le comptage des lignes de flots négatives du gradient de $f$. Bien que le complexe lui-même dépende fortement de la fonction de Morse choisie, l'homologie du complexe est
indépendante de ce choix. Le papier sera organisé comme suit: dans l'optique de calculer le complexe de chaîne défini sur les points critiques d'une fonction de Morse $f$ et l'homologie correspondante dans
le cas des coefficients entiers\footnote{Il est à noter que le signe associé à une orientation est, en pratique, une problématique bien subtile lorsque nous calculons l'opérateur de bord induisant l'homologie de Morse. Afin d'éviter de tels problèmes techniques, l'homologie est généralement calculée sur $\mathds{Z}/2\mathds{Z}$ au lieu de $\mathds{Z}$. Or, c'est le dernier cas qui sera considéré dans le présent article.}, nous définirons, en premier lieu, correctement les orientations
des espaces de modules des flots et nous fournirons une preuve de la cohérence de cette orientation. Enfin, nous appliquerons cette construction et l’utiliserons lors du calcul explicite de ces homologies sur une variété de
dimension 4 -- i.e. $S^2\times T^2$. Nous conclurons la discussion avec une invitation à généraliser la discussion dans un contexte de dimension infinie. \\

Avant de clore cette section d'introduction, il est également pertinent de mentionner la grande utilité de la théorie Morse (ainsi que sa version complexe, la théorie de Picard-Lefschetz \cite{vassiliev2002applied}) dans le contexte de la physique théorique \cite{daemi2019equivariant, Frellesvig_2019, witten2010new}. En effet, il s'avère que la théorie Morse a une formulation physique très intéressante; par exemple, une fonction Morse peut être considérée comme une sorte de potentiel, donc le flux induit par son gradient est la force subie par une particule (ou une corde). Les points d'équilibre de la particule correspondent alors aux points critiques du potentiel. Si nous prenons une extension supersymétrique de ce potentiel \cite{wess1992supersymmetry}, interprétant le tout comme une théorie quantique, la structure de l'état fondamental calcule l'homologie de l'espace dans lequel les particules (ou les cordes) se déplacent. Ceci reflète bien le fait qu'il est possible de déduire de l'information sur la topologie d'une variété en étudiant la dynamique d'une (super)corde soumise à un (super)potentiel -- voir \cite{freed1999five} pour une couverture complète et plusieurs exemples. La théorie de Morse apparait aussi dans de récents développements sur les intégrales de Feynman et les intégrales de cordes -- e.g. \cite{mizera2019aspects,mizera2020status} -- et la topologie de cordes \cite{cohen2004morse}. \\
\section{Remerciements}
L'auteur souhaite remercier Pr. François Lalonde pour lui avoir donné l'opportunité de faire cette étude, ainsi que pour ses nombreux conseils et astuces sans lesquels ce texte ne se serait pas écrit. L'auteur est aussi reconnaissant envers l'Université de Montréal pour son hospitalité lors de l'été 2018.
\section{Orientation et cohérence}
Dans ce qui suivra, nous assumerons que le couple $(f,X)$ satisfait la condition de Morse-Smale \cite{audin2014morse,banyaga2013lectures}. Il clair que pour tout point critique $\xi$ d'une variété orientée $M$ avec $\mathcal{O}_M(\xi)$ la variété stable $W^s(\xi)$ et la variété instable $W^u(\xi)$ sont difféomorphiques à un $p$-disque pour $0\le p\le \dim(M)$. Ainsi, ces sous-variétés de $M$ sont orientables. La convention adoptée ci-bas sera la suivante. L'orientation $\mathcal{O}_{W^u}(\xi)$ sera définie en premier pour tous les points critiques $\xi$ tels que 
\begin{equation}
\mathcal{O}_{W^u}(\xi)\oplus \mathcal{O}_{W^s}(\xi)=\mathcal{O}_{M}(\xi).
\end{equation}
Par conséquent, une orientation définie sur $W^u(\xi)$ en induit une naturelle sur $W^s(\xi)$.\\

Voici la proposition principale sur l'orientation des modules de flots.
\begin{prop} \textnormal{(Orientations induites)} 
$ \forall \ \xi\in \textnormal{Crit}(f)\ : \textnormal{Ind}(\xi)>0$ fixons une orientation arbitraire sur $W^u(\xi)$ que l'on note $\mathcal{O}_{W^u}(\xi)$. Ainsi, $\forall \ \xi, \gamma\in \textnormal{Crit}(f)$ la variété connectée 
\begin{equation}
_\xi M_\gamma := W^u(\xi)\cap W^s(\gamma),
\end{equation}
et l'espace des orbites associé\footnote{Une définition équivalente, mais sans doute plus explicite, de $_\xi\mathcal{M}_{\gamma}$ serait $$_\xi\mathcal{M}_\gamma:=\{ \textnormal{ligne de flots de} \ \Grad f \ \textnormal{reliant} \ \xi \ \textnormal{et} \ \gamma \}/\textnormal{reparamétrisation},$$
où le quotient par \emph{reparamétrisation} identifie $u(s)$ à $u(s+\textnormal{constante})$ et où le (pseudo)gradient, dénoté $\Grad$, est introduit plus bas dans le texte.} 
\begin{equation}
_\xi \mathcal{M}_\gamma := W^u(\xi)\cap W^s(\gamma)\cap f^{-1}(\eta),
\end{equation}
où $\eta$ est une valeur régulière de $f$ dans $(f(\xi),f(\gamma))$, disposent d'orientations induites, respectivement, $\mathcal{O}_{_\xi M_\gamma }^{\textnormal{ind}}$ et $\mathcal{O}_{_\xi \mathcal{M}_\gamma }^{\textnormal{ind}}$.
\end{prop}
Nous prouverons les lemmes suivants en préalable à la preuve de la proposition centrale.
\begin{lemma}
Soit $M$ une variété orientée. Soient $S_1$ et $S_2$ deux sous-variétés orientées de $M$ transverses. Alors, la sous-variété $S_1\cap S_2$ est orientée.
\end{lemma}
\begin{prL}
Ceci est clair puisque la transversalité de l'intersection -- i.e. $S_1\pitchfork S_2$ -- implique que
\begin{equation}
\forall \ z\in S_1\cap S_2, \mathcal{T}_z(S_1\cap S_2)=\mathcal{T}_z S_1\cap \mathcal{T}_z S_2,
\end{equation}
(preuve complète en annexe) et de l'algorithme de Zassenhaus \cite{luks1997some} utilisant les orientations sur les espaces tangents respectifs. Q.E.D
\end{prL}
Comme $ W^u(\xi)$ et $ W^s(\gamma)$ sont des sous-variétés orientées et dont l'intersection est transversale sous la condition que $X$ est de Smale, il s'en suit que $W^u(\xi)\cap W^s(\gamma)$ est aussi orientée.




\begin{lemma} \textnormal{(Caractérisation de $\mathcal{O}$ des sous-variétés d'une variété orientée)} 
Soit $M$ une variété orientée. Soit $S$ une $n$-sous-variété de codimension $k$ dans $M$. Alors, $S$ est orientée si et seulement si son fibré normal \cite{snaith},
que nous nommerons $NS$, est un fibré vectoriel orienté. 
\end{lemma}
\begin{prL} $S$ est de dimension $n$ et de codimension $k$.
Restreint à $S$ (aux points de $S$ dans $M$), le fibré tangent de $M$ se décompose comme suit 
\begin{equation}
TM|_{S}\cong TS\oplus NS.
\end{equation}
De l'algèbre multilinéaire, nous avons l'identité\footnote{Soient $V$ et $W$ deux espaces vectoriels sur un corps $\mathds{K}$. Alors, rappelons que 
$$\bigwedge(V\oplus W)\cong\bigwedge(V)\otimes_{\mathds{K}}\bigwedge(W).$$}
\begin{equation}
\bigwedge^{n+k}(TM|_S)=\bigwedge^{n+k}(TS\oplus NS)\cong\bigwedge^n TS \otimes_{\mathds{R}}\bigwedge^k NS,
\end{equation}
où $n$ est la puissance extérieure maximale du fibré $TS$ et $k$ puissance extérieure maximale du fibré $NS$. Comme, par hypothèse, $M$ est orientée, par définition d'un fibré vectoriel orienté, la top-wedge puissance $\bigwedge^{n+k}(TM)$ est un fibré en droites trivialisable. Il est clair que la restriction $\bigwedge^{n+k}(TM|_S)$ l'est aussi. Ainsi, $TM|_S$ est orienté.\\
"$\Rightarrow$" Assumons que $S$ est une sous-variété orientée. En suivant l'idée de la preuve du premier lemme, il est clair que $TS$ est aussi orienté. Par la relation en somme directe, on induit une orientation sur $NS$.\\
"$\Leftarrow$" Assumons que $NS$ est un fibré vectoriel orienté. Par la relation en somme directe, comme $M$, donc $TM|_S$, est orientée par hypothèse, cela induit une orientation sur $TS$ et donc sur $S$. Q.E.D
\end{prL}
\begin{lemma} \textnormal{(Lemme de séparation)} donné une chaîne courte et exacte munie des applications $i$ et $\pi$ entre objets de catégories 
\begin{equation}
0 \rightarrow A \xrightarrow{i} B \xrightarrow{\pi} C \rightarrow 0,
\end{equation}
alors $B\cong A\oplus C$ avec $i$ défini comme l'injection canonique de $A$ et $\pi$ la projection canonique sur $C$ -- i.e. il existe un isomorphisme de chaînes courtes et exactes  
\begin{equation}
(0 \rightarrow A \xrightarrow{i} B \xrightarrow{\pi} C \rightarrow 0)\cong (0 \rightarrow A \xrightarrow{i} A\oplus C \xrightarrow{\pi} C \rightarrow 0).
\end{equation}
\end{lemma}
\begin{prL}
Voir \cite{weibel1995introduction}.
\end{prL}
\begin{lemma}
Soient $A,B,C$ des $R$-modules. Soit $C=A+B$. Alors, $A\cong C/B$.
\end{lemma}
\begin{prL}
Établissons la chaîne courte suivante par le truchement de l'inclusion $((i: B\hookrightarrow A+B):\{y\mid y\in B\}\to \{y+0\mid y \in B, 0\in A\})$ (injective) et la projection $((\pi:A+B\to A):\{x+y\mid x\in A, y\in B\} \to \{(x+y)-y\mid x\in A\})$ (surjective) 
\begin{equation}
0 \rightarrow B \xrightarrow{i} A+B \xrightarrow{\pi} A \rightarrow 0.
\end{equation}
Par construction, nous remarquons : (i) $i(x)=x \Rightarrow \textnormal{Im}(i)=B$, (ii) $\ker(\pi)=B$. Ainsi, $\ker(\partial_k)\equiv \textnormal{Im}(\partial_{k+1})$. Il s'en suit que cette chaîne est courte et exacte. De plus, comme $\pi$ est surjective, par le premier théorème d'isomorphisme, 
\begin{equation}
\begin{split}
A&\cong \frac{A+B}{\ker(\pi)}=\frac{A+B}{\textnormal{Im}(i)}\cong \frac{A+B}{B}=\frac{C}{B}.
\end{split}
\end{equation}
Q.E.D
\end{prL}
Dans ce qui suivra $R=\mathds{R}$, comme $\mathds{R}$ est un corps et donc un anneau, le $R$-module équivaut à un $R=\mathds{R}$-espace vectoriel. Par fibre, la condition de Smale sur $X$ implique une relation de transversalité qui se réduit, par abus de notation sur la restriction de ces fibrés sur $_\xi M_\gamma$ \footnote{La condition de transversalité est satisfaite pour tout point dans $_\xi M_\gamma$, sous-variété de $M$.} et à l'aide du lemme ci-haut, à
\begin{equation}
TW^u(\xi)|_{_\xi M_\gamma}+TW^s(\gamma)|_{_\xi M_\gamma}= T_\xi M_\gamma \Rightarrow TW^u(\xi)|_{_\xi M_\gamma}\cong T_\xi M_\gamma/TW^s(\gamma)|_{_\xi M_\gamma}\cong NW^s(\gamma)|_{_\xi M_\gamma}.
\end{equation}
\begin{lemma} \textnormal{(Relation en somme directe)} L'isomorphisme 
\begin{equation}
TW^u(\xi)|_{_\xi M_\gamma}\cong T _\xi M_\gamma\oplus NW^s(\gamma)|_{_\xi M_\gamma},
\end{equation}
existe. 
\end{lemma}
Avant de donner la preuve, rappelons que comme $_\xi M_\gamma:=W^u(\xi)\cap W^s(\gamma)$ est transverse, nous avons, par le lemme en annexe,
\begin{equation}
T _\xi M_\gamma:=T(W^u(\xi)\cap W^s(\gamma))=TW^u(\xi)\cap TW^s(\gamma).
\end{equation}
Ainsi, il existe l'application injective canonique 
\begin{equation}
(i:T _\xi M_\gamma=TW^u(\xi)\cap TW^s(\gamma)\to \tilde{T}W^u(\xi)):\Theta \to \Theta|_\cap,
\end{equation}
où $\tilde{T}W^u(\xi)=TW^u(\xi)|_{TW^u(\xi)\cap TW^s(\gamma)}$. Muni de cette injection nous pouvons fournir la preuve suivante.
\begin{prL}
Prenons la chaîne courte suivante 
\begin{equation}
0 \rightarrow T _\xi M_\gamma \xrightarrow{i} TW^u(\xi)|_{_\xi M_\gamma} \xrightarrow{\pi} NW^s(\gamma)|_{_\xi M_\gamma} \rightarrow 0,
\end{equation}
où $i$ est l'injection ci-haute et la projection $\pi$, ici l'identité, puisque nous avons une équivalence entre chaînes 
\begin{equation}
0 \rightarrow T _\xi M_\gamma \xrightarrow{i} TW^u(\xi)|_{_\xi M_\gamma} \xrightarrow{\pi} NW^s(\gamma)|_{_\xi M_\gamma} \rightarrow 0,
\end{equation}
et
\begin{equation}
0 \rightarrow T _\xi M_\gamma \xrightarrow{i} TW^u(\xi)|_{_\xi M_\gamma} \xrightarrow{\pi} TW^u(\xi)|_{_\xi M_\gamma} \rightarrow 0,
\end{equation}
dû au fait que $X$ soit de Smale -- i.e. $TW^u(\xi)|_{_\xi M_\gamma}\cong NW^s(\gamma)|_{_\xi M_\gamma}$. Par conséquent, $\ker(\pi)=0$. Cela implique que les groupes d'homologie soient triviaux et que la chaîne soit exacte. Du lemme de séparation, nous avons le résultat désiré. Q.E.D
\end{prL}
\begin{prP}
Clairement, du lemme 1, comme $X$ est de Smale, $W^u(\xi)\cap W^s(\gamma)$ est orientable. Ainsi $_\xi M_\gamma$ est orientable. Dans ce qui suit, nous observerons comment induire une orientation sur la variété connectée $_\xi M_\gamma$ en obtenant d'abord une orientation pour le fibré tangent de la variété instable de $\xi$ et une orientions pour le fibré normal de la variété stable de $\gamma$ en imposant une restriction de ces fibrés par la variété connectée elle-même. De façon analogue, nous obtiendrons une orientation pour l'espace des orbites correspondant via l'orientation induite plus tôt sur la variété connectée et l'orientation naturelle du fibré en droites de $\mathds{R}$, donné le pseudo-gradient de $f$.\\

Du travail fait plus haut, 
\begin{equation}
\label{eq1}
TW^u(\xi)|_{_\xi M_\gamma}\cong T _\xi M_\gamma \oplus NW^s(\gamma)|_{_\xi M_\gamma}.
\end{equation}
Le but est maintenant de montrer que $TW^u(\xi)|_{_\xi M_\gamma}$ et que $NW^s(\gamma)|_{_\xi M_\gamma}$ sont deux fibrés orientés, rendant ainsi évidant une orientation sur $T _\xi M_\gamma$ et donc sur $_\xi M_\gamma$.\\

En premier lieu, la restriction du fibré total $TW^u(\xi)$ à n'importe quelle sous-variété de $M$ -- e.g. $_\xi M_\gamma$ -- doit être orientée. Ainsi, $TW^u(\xi)|_{_\xi M_\gamma}$ est un fibré orienté. En second lieu, une orientation sur $W^u(\gamma)$ induit par construction une orientation sur $W^s(\gamma)$. Ceci tient si et seulement si $NW^s(\gamma)$, donc $NW^s(\gamma)|_{_\xi M_\gamma}$, est un fibré orienté. Explicitement, on restreint le fibré vectoriel 
\begin{equation}
\pi : NW^s(\gamma)\to W^s(\gamma),
\end{equation}
à (une fibre de) la sous-variété $_\xi M_\gamma \hookrightarrow W^s(\gamma)$, pour obtenir une orientation sur le fibré normal restreint à la variété connectée. Par la relation en somme directe \eqref{eq1}, ces deux orientations en induisent une sur le fibré tangent de la variété connectée de $M$ entre $\xi$ et $\gamma$ et, par définition, sur $_\xi M_\gamma$.\\

Nous voulons maintenant induire une orientation sur l'espace des orbites associé $_\xi \mathcal{M}_\gamma$. Il y a une décomposition en somme directe très évidente reliant le fibré tangent de cet espace au fibré tangent de la variété connectée de $M$ entre $\xi$ et $\gamma$ -- i.e. 
\begin{equation}
T _\xi M_\gamma|_{_\xi \mathcal{M}_\gamma}\cong \mathds{R}\oplus T _\xi \mathcal{M}_\gamma.
\end{equation}
Notons que de plus haut, la restriction $T _\xi M_\gamma|_{_\xi \mathcal{M}_\gamma}$ est orientée. De plus, donné le pseudo-gradient, dénoté $\Grad(f)$ et défini comme dans \cite{audin2014morse},
$\mathds{R}$ est clairement un fibré en droite orienté. De la dernière relation en somme directe, $T _\xi \mathcal{M}_\gamma$ est donc orienté et $ _\xi \mathcal{M}_\gamma$ par le fait même. Q.E.D
\end{prP}
\textbf{Remarque}: Comme le note \cite[p.~39]{audin2014morse} $\mathcal{L}(\xi,\gamma)$, l'espace des modules de flots entre $\xi$ et $\gamma$, correspond à $_\xi \mathcal{M}_\gamma$. Ainsi, on induit une orientation sur $\mathcal{L}(\xi,\gamma)$.\\ 

Nous prouvons maintenant la proposition suivante. 
\begin{prop}\label{prop2} \textnormal{(Cohérence)} L'application de collage $\mathcal{O}^{\omega_\eta}$, induite sur les orientations des fibrés par $\mathfrak{S}$ (définie plus bas) avec $\omega_\eta:=u\mathfrak{S}_\eta v$ ainsi que les orientations obtenues dans la première proposition sont compatibles -- i.e. 
\begin{equation}
\mathcal{O}^{\omega_\eta}(\mathcal{O}(_\xi M_\gamma ^u)_{\textnormal{Ind}},\mathcal{O}(_\gamma M_\zeta ^v)_{\textnormal{Ind}})=\mathcal{O}(_\xi M_\zeta ^{\omega_\eta}).
\end{equation}
\end{prop}
Avant d'élaborer une preuve, nous avons besoin des préliminaires suivants. 
\begin{lemma} \textnormal{(Application de collage)} Soit $(f,X)$ un couple de Morse-Smale. Soient $\xi,\gamma, \zeta\in \textnormal{Crit}(f)$ tels que $\textnormal{Ind}(\xi)=k+1$, $\textnormal{Ind}(\gamma)=k$ et $\textnormal{Ind}(\zeta)=k-1$. Alors, il existe un nombre réel positif $\eta_0$ et un plongement (immersion injective)
\begin{equation}
\mathfrak{S}: _\xi \mathcal{M}_\gamma \times [\eta_0,\infty[\times _\gamma \mathcal{M}_\zeta \to _\xi \mathcal{M}_\zeta,
\end{equation}
tel que 
\begin{equation}
(u,\eta,v)\mapsto u\mathfrak{S}_\eta v:=\omega_\eta,
\end{equation}
et 
\begin{equation}
\omega_\eta \xrightarrow{\eta \to \infty} (u,v) \ \ \ \ \textnormal{et} \ \ \ \ \omega_\eta \xrightarrow{\eta\to \eta_0}(u,v)|_{\textnormal{passe par} \ \gamma}.
\end{equation}
Notons que la dernière condition est équivalente, donnée la métrique de l'espace ambiant $d$, à $d(\omega_\eta,\gamma)\to 0$.
\end{lemma}
\begin{prL}
Voir \cite[Proposition 2.56]{schwarz1993morse}. Q.E.D
\end{prL}
Avec les éléments des annexes en main, nous sommes prêts à considérer la construction suivante.\\

donnés $\xi,\gamma,\zeta \in \textnormal{Crit}(f)$ d'indices de Morse $k+1$, $k$ et $k-1$ respectivement et donnés $u\in  _\xi \mathcal{M}_\gamma$ et $v\in _\gamma \mathcal{M}_\zeta$, nous verrons que l'application de collage $\mathfrak{S}$ induit une application de collage des orientations -- i.e. 
\begin{equation}
\mathcal{O}^{\omega_\eta}: \textnormal{Or}(_\xi M_\gamma ^u)\times \textnormal{Or}(_\gamma M_\zeta ^v)\to \textnormal{Or}(_\xi M_\zeta ^{\omega_\eta}),
\end{equation}
où $_\xi M_\gamma ^u$ dénote la composante connectée de $\xi M_\gamma$ contenant $u$ (voir la figure ci-haut). Par construction, la courbe décrite par le flot $\phi_\eta(u)$ -- i.e. $_\xi M_\gamma ^u$ de dimension 1 -- satisfait $\frac{\textnormal{d}}{\textnormal{d}\eta}\phi_\eta(u)\neq 0$ comme aucun point critique n'est atteint. Ainsi, comme $\frac{\textnormal{d}}{\textnormal{d}\eta}\phi_\eta(u)=\Grad(f(\phi_\eta(u)))\neq 0$, une orientation sur $_\xi M_\gamma^u$ est alors induite par le pseudo gradient, dénotée par $\mathcal{O}^{\dot{\phi}_\eta}(u)$,
et de façon similaire pour $_\gamma M_\zeta ^v$, par
$\mathcal{O}^{\dot{\phi}_\eta}(v)$. De la façon avec laquelle nous avons construit la preuve sur le plongement $\mathfrak{S}$, nous avons que les vecteurs des champs de vecteurs $\Grad f(\omega_\eta)$ et $\frac{\textnormal{d}}{\textnormal{d}\eta}\omega_\eta$ sont linéairement indépendants. Souvenons-nous que $\dim(_\xi M^{\omega_\eta}_\zeta)=2$. Pour chaque $\eta\in[\eta_0,\infty[$, $(\Grad f(\omega_\eta), -\frac{\textnormal{d}}{\textnormal{d}\eta}\omega_\eta)$ est un cadre de champs vectoriels \cite{loring2011introduction} qui génère l'espace bidimensionnel $_\xi M^{\omega_\eta}_\zeta$. 
Ainsi, $(\Grad f(\omega_\eta), -\frac{\textnormal{d}}{\textnormal{d}\eta}\omega_\eta )$ forme un cadre global continu 
donnant lieu a à une orientation sur l'espace $_\xi M_\zeta ^{\omega_\eta}$. Cette orientation sera dénotée 
\begin{equation}
\mathcal{O}(\Grad f(\omega_\eta), -\frac{\textnormal{d}}{\textnormal{d}\eta}\omega_\eta).
\end{equation}
De façon explicite, nous définissons alors l'application de collage de la proposition~\ref{prop2} comme suit
\begin{equation}
\mathcal{O}^{\omega_\eta}(\mathcal{O}^{\dot{\phi}_\eta}(u),\mathcal{O}^{\dot{\phi}_\eta}(v)):=\mathcal{O}(\Grad f(\omega_\eta), -\frac{\textnormal{d}}{\textnormal{d}\eta}\omega_\eta).
\end{equation}
De la géométrie différentielle, nous avons qu'une variété lisse, connectée et orientable a exactement deux orientations. Ainsi, pour un cas général et non basé sur l'orientation induite par le flot, il semble approprié de définir l'application de collage comme 
\begin{equation}
\begin{split}
\mathcal{O}^{\omega_\eta}(\mathcal{O}(_\xi M_\gamma ^u),\mathcal{O}(_\gamma M_\zeta ^v))&:=\alpha \beta \mathcal{O}^{\omega_\eta}(\mathcal{O}^{\dot{\phi}_\eta}(u),\mathcal{O}^{\dot{\phi}_\eta}(v))\\&=\alpha\beta \mathcal{O}(\Grad f(\omega_\eta), -\frac{\textnormal{d}}{\textnormal{d}\eta}\omega_\eta),
\end{split}
\end{equation}
où $\alpha,\beta\in\{\pm 1\}$ peuvent être déterminés par le truchement des relations suivantes 
\begin{equation}
\mathcal{O}(_\xi M_\gamma ^u)=\alpha \mathcal{O}^{\dot{\phi}_\eta}(u) \ \ \ \textnormal{et} \ \ \ \mathcal{O}(_\gamma M_\zeta ^v)=\beta \mathcal{O}^{\dot{\phi}_\eta}(v).
\end{equation}
Avec ceci, nous pouvons attaquer la preuve de la seconde proposition. 
\begin{prP}
Définissons $n_u\in \{\pm 1\}$ par la relation $\mathcal{O}(_\xi M_\gamma ^u)_{\textnormal{Ind}}=n_u \mathcal{O}^{\phi_\eta}(u)$ et similairement pour $v$. Alors, par définition de l'application de collage $\mathcal{O}^{\omega_\eta}$
\begin{equation}
\begin{split}
\mathcal{O}^{\omega_\eta}(\mathcal{O}(_\xi M_\gamma ^u)_{\textnormal{Ind}},\mathcal{O}(_\gamma M_\zeta ^v)_{\textnormal{Ind}})&=n_un_v\mathcal{O}^{\omega_\eta}(\mathcal{O}^{\dot{\phi}_\eta}(u),\mathcal{O}^{\dot{\phi}_\eta}(v))\\&=n_un_v\underbrace{\mathcal{O}(\Grad f(\omega_\eta), -\frac{\textnormal{d}}{\textnormal{d}\eta}\omega_\eta)}_{(\star)}.
\end{split}
\end{equation}
Notre but est évidemment de construire une façon de comparer $(\star)$ avec une orientation induite sur $_\xi M_\zeta ^{\omega_\eta}$, $\mathcal{O}(_\xi M_\zeta ^{\omega_\eta})_{\textnormal{Ind}}$ --- obtenue comme suit.\\

La première chose à faire, ce sera de relier les orientations induites sur les fibrés tangents $T _\xi M_\gamma ^u$, $T _\gamma M_\zeta ^v$ et $T _\xi M_\zeta ^{\omega_\eta}$. Malheureusement, comme $W^s(\gamma)\cap W^u(\gamma)=\emptyset$ les bases des fibrés -- i.e. $_\xi M_\gamma ^u$, $_\gamma M_\zeta ^v$ et $_\xi M_\zeta ^{\omega_\eta}$ -- n'ont pas nécessairement un point commun. Or, par l'existence des limites \ref{appC} $\mathfrak{E}^\pm$, il est clair que le singleton $\{\gamma\}$ réside dans la frontière de toutes ces bases. Ainsi, les bases ont un point commun si nous étendons les fibré à $\{\gamma\}$ en considérant la frontière de celles-ci. Ces extensions seront notées $T^\gamma _\xi M_\gamma ^u$, $T^\gamma _\gamma M_\zeta ^v$ et $T^\gamma _\xi M_\zeta ^{\omega_\eta}$ et \emph{seront nécessaires} pour montrer comment $_\xi M_\zeta ^{\omega_\eta}$ est orienté. Rappelons de \ref{appC} que les orientations correspondantes à ces limites sont dénotées $\mathcal{O}(\mathfrak{E}^\pm)$.\\

De la première proposition 
\begin{equation}
T W^u(\xi)|_{_\xi M_\gamma ^u} \cong T _\xi M_\gamma ^u \oplus N W^s(\gamma)|_{_\xi M_\gamma ^u}. 
\end{equation}
Par extension sur $\{\gamma\}$, cet isomorphisme en somme directe se réécrit, comme les extensions de $W^u(\xi)$ et $W^s(\gamma)$ à $\{\gamma\}$ est l'identité puisque $\gamma$ en est un élément \footnote{Ici, aucune adaptation de la notation n’est nécessaire pour le premier et le dernier terme de l'équation \eqref{refran} -- or, nous ajouterons $\gamma$ en superscript à $T\star$ et $N\star$ afin de garder une consistance à la notation utilisée pour l'extension des fibrés à $\{\gamma\}$.}, 
\begin{equation}
\label{refran}
T^\gamma W^u(\xi)|_{_\xi M_\gamma ^u} \cong T^\gamma _\xi M_\gamma ^u \oplus N^\gamma W^s(\gamma)|_{_\xi M_\gamma ^u} ,
\end{equation}
générant la relation en somme directe sur les orientations
\begin{equation}
\mathcal{O}(T^\gamma W^u(\xi)|_{_\xi M_\gamma ^u}) = \mathcal{O}(T^\gamma _\xi M_\gamma ^u)_{\textnormal{Ind}} \oplus \mathcal{O}(N^\gamma W^s(\gamma)|_{_\xi M_\gamma ^u}).
\end{equation}
En restreignant l'égalité à la fibre associée à $\gamma$
\begin{equation}
\label{eq2}
\mathcal{O}(T^\gamma W^u(\xi)|_{\gamma}) = \mathcal{O}(T^\gamma _\xi M_\gamma ^u|_\gamma)_{\textnormal{Ind}} \oplus \mathcal{O}(N^\gamma W^s(\gamma)|_{\gamma }).
\end{equation}
De façon analogue, il est vrai de la première proposition que 
\begin{equation}
T W^u(\gamma)|_{_\gamma M_\zeta ^v} \cong T _\gamma M_\zeta ^v \oplus N W^s(\zeta)|_{_\gamma M_\zeta ^v},
\end{equation}
et ainsi 
\begin{equation}
\mathcal{O}(T W^u(\gamma)|_{_\gamma M_\zeta ^v}) = \mathcal{O}(T _\gamma M_\zeta ^v)_{\textnormal{Ind}} \oplus \mathcal{O}(N W^s(\zeta)|_{_\gamma M_\zeta ^v}).
\end{equation}
En considérant l'extension\footnote{Ici, l'adaptation de notation n'est plus abusive, dans le sens où l'extension ne correspond plus nécessairement à l'identité.} à $\{\gamma\}$
\begin{equation}
\mathcal{O}(T^\gamma W^u(\gamma)|_{_\gamma M_\zeta ^v}) = \mathcal{O}(T^\gamma _\gamma M_\zeta ^v)_{\textnormal{Ind}} \oplus \mathcal{O}(N^\gamma W^s(\zeta)|_{_\gamma M_\zeta ^v}).
\end{equation}
Ayant l'extension jusqu’à $\{\gamma\}$ en main, il est possible de restreindre l'égalité à la fibre correspondante à $\gamma$
\begin{equation}
\label{eq3}
\mathcal{O}(T^\gamma W^u(\gamma)|_{\gamma}) = \mathcal{O}(T^\gamma _\gamma M_\zeta^v|\gamma)_{\textnormal{Ind}} \oplus \mathcal{O}(N^\gamma W^s(\zeta)|_{\gamma}).
\end{equation}
De plus, remarquons que par définition
\begin{equation}
W^s(\gamma)\cap W^u(\gamma)=\emptyset \Rightarrow W^s(\gamma)\pitchfork W^u(\gamma),
\end{equation}
et alors 
\begin{equation}
TW^s(\gamma)\cap TW^u(\gamma)=\emptyset. 
\end{equation}
De cette façon, sur cette intersection $"+ \cong \oplus"$. Par hypothèse de transversalité 
\begin{equation}
TW^u(\gamma)+ TW^u(\gamma)= TM \Leftrightarrow TW^u(\gamma)\oplus TW^u(\gamma)\cong TM.
\end{equation}
Par définition du fibré normal, 
\begin{equation}
TM\cong NW^s(\gamma) \oplus TW^s(\gamma).
\end{equation}
De ces deux dernières équations, nous obtenons l'isomorphisme suivant 
\begin{equation}
TW^u(\gamma)\cong N W^s(\gamma)\Rightarrow T^\gamma W^u(\gamma)|_{\gamma}\cong N^\gamma W^s(\gamma)|_{\gamma}.
\end{equation}
En utilisant \eqref{eq2}, \eqref{eq3} et le dernier isomorphisme, il est possible d'écrire
\begin{equation}
\begin{split}
\mathcal{O}(T^\gamma W^u(\xi)|_{\gamma})& = \mathcal{O}(T^\gamma _\xi M_\gamma ^u|_\gamma)_{\textnormal{Ind}} \oplus \mathcal{O}(N^\gamma W^s(\gamma)|_{\gamma})\\&=\mathcal{O}(T^\gamma _\xi M_\gamma ^u|_\gamma)_{\textnormal{Ind}} \oplus \mathcal{O}(T^\gamma W^u(\gamma)|_{\gamma})\\&=\mathcal{O}(T^\gamma _\xi M_\gamma ^u|_\gamma)_{\textnormal{Ind}} \oplus \mathcal{O}(T^\gamma _\gamma M_\zeta ^v|_\gamma)_{\textnormal{Ind}}\oplus \mathcal{O}(N^\gamma W^s(\zeta)|_{\gamma}),
\end{split}
\end{equation}
Comme $\omega_\eta \xrightarrow{\eta\to \eta_0}(u,v)|_{\textnormal{passe par} \ \gamma}$, dans cette limite 
\begin{equation}
\begin{split}
\mathcal{O}(T^\gamma W^u(\xi)|_{\gamma})& = \mathcal{O}(T^\gamma _\xi M_\zeta ^{\omega_\eta}|_\gamma)_{\textnormal{Ind}}\oplus \mathcal{O}(N^\gamma W^s(\zeta)|_{\gamma}).
\end{split}
\end{equation}
Ceci prouve que $T^\gamma _\xi M_\zeta ^{\omega_\eta}|_\gamma$ est une variété orientable et en définit même l'une des deux orientations. De la discussion dans \cite[p.~104]{hirsch2012differential} comme $_\xi M_\zeta ^{\omega_\eta}$ est connectée nous avons qu'une orientation sur le fibré total -- i.e. $T^\gamma _\xi M_\zeta ^{\omega_\eta}$ -- est déterminée par l'orientation d'une seule fibre arbitraire 
-- i.e. ici $T^\gamma _\xi M_\zeta ^{\omega_\eta}|_\gamma$. Ainsi, $_\xi M_\zeta ^{\omega_\eta}$ est une variété orientable.\\ 

Cela dit, nous avons la décomposition naturelle de l'orientation en somme directe suivante en extensions de $\{\gamma\}$
\begin{equation}
\mathcal{O}(T^\gamma _\xi M_\zeta ^{\omega_\eta})_{\textnormal{Ind}}=\mathcal{O}(T^\gamma _\xi M_\gamma ^u)_{\textnormal{Ind}}\oplus \mathcal{O}(T^\gamma _\gamma M_\zeta ^v)_{\textnormal{Ind}}= n_u\mathcal{O}(\mathfrak{E}^+)\oplus n_v \mathcal{O}(\mathfrak{E^-}).
\end{equation}
Comme une variété orientable possède deux orientations possibles, nous devons avoir la correspondance suivante 
\begin{equation}
n_u\mathcal{O}(\mathfrak{E}^+)\oplus n_v \mathcal{O}(\mathfrak{E^-})=n_un_v\mathcal{O}(\Grad f(\omega_\eta), -\frac{\textnormal{d}}{\textnormal{d}\eta}\omega_\eta).
\end{equation}
Cette dernière condition sur l'orientation des fibrés implique 
\begin{equation}
\mathcal{O}^{\omega_\eta}(\mathcal{O}(_\xi M_\gamma ^u)_{\textnormal{Ind}},\mathcal{O}(_\gamma M_\zeta ^v)_{\textnormal{Ind}})=\mathcal{O}(_\xi M_\zeta ^{\omega_\eta}).
\end{equation}
Q.E.D
\end{prP}
\section{Calcul des groupes d'homologie de Morse pour $S^2\times T^2$}
Dans le contexte de la théorie de Morse, nous désirons calculer les groupes d'homologies du produit $S^2\times T^2$.\\

De façon générale, soient $M_1$ et $M_2$ deux variétés munies de fonctions de Morse $f_1$ et $f_2$, respectivement, et munies de champs pseudo-gradients $X_1$ et $X_2$ satisfaisant la condition de Smale. Il est à noter que 
\begin{equation}
f_1+f_2: M_1\times M_2 \to \mathds{R},
\end{equation}
doit être une fonction de Morse. Les points critiques de $f_1+f_2$ sont les points $(\xi_1,\xi_2)$, où $\xi_1\in \textnormal{Crit}(f_1)$ et $\xi_2\in \textnormal{Crit}(f_2)$. Supposons le contraire. Alors, comme la différentielle $f_{i_\star}:\mathcal{T}_p M_i\to \mathds{R}$ sur $f_i$ au point $p\in M_i$ est linéaire, 
\begin{gather}
(f_1+f_2)_\star=f_{1_\star}+f_{2_\star}=0 \ \lightning, \\ \therefore \ (\xi_1,\xi_2)\in \textnormal{Crit}_k(f_1+f_2)\Leftrightarrow (\xi_1\in \textnormal{Crit}_i(f_1) \ \wedge \ \xi_2\in \textnormal{Crit}_j(f_2):i+j=k).
\end{gather}
Hormis cela, comme $M_1$ et $M_2$ sont deux espaces différents, il est clair que le couple $(X_1,X_2)$ est lui aussi de Smale. Soit alors un point critique $\xi:=(\xi_1,\xi_2)$ de $f_1+f_2$ d'indice $k$ et soit un point critique $\gamma:=(\gamma_1,\gamma_2)$ de la même fonction, mais cette fois-ci d'indice $k-1$ tel qu'il existe \emph{une} trajectoire connectant $\xi$ à $\gamma$. Ce flot est décrit par 
\begin{equation}
\phi_t^{(X_1,X_2)}(\cdot, \star):=(\phi_t^{X_1}(\cdot),\phi_t^{X_2}(\star)),
\end{equation}
de sorte que 
\begin{equation}
(_\xi \mathcal{M}_\gamma)_{(X_1,X_2)}\cong (_{\xi_1} \mathcal{M}_{\gamma_1})_{X_1}\times  (_{\xi_2} \mathcal{M}_{\gamma_2})_{X_2}.
\end{equation}
Remarquons que si $\xi_1\neq \gamma_1$ et si $\xi_2\neq \gamma_2$, pour avoir $(_\xi \mathcal{M}_\gamma)_{(X_1,X_2)}\neq \emptyset$, nous demandons 
\begin{equation}
(\textnormal{Ind}(\xi_1)\ge \textnormal{Ind}(\gamma_1)+1) \ \wedge \ (\textnormal{Ind}(\xi_2)\ge \textnormal{Ind}(\gamma_2)+1).
\end{equation}
Ainsi, en général, ces conditions se contractent comme
\begin{equation}
\textnormal{Ind}(\xi)\ge \textnormal{Ind}(\gamma)+2.
\end{equation}
Si, dans un autre contexte, $\xi$ et $\gamma$ sont d'indices consécutifs (alors soit que $\xi_1=\gamma_1$ ou soit que $\xi_2=\gamma_2$), nous avons 
\begin{equation}
(_\xi \mathcal{M}_\gamma)_{(X_1,X_2)}\cong \begin{cases} 
      \xi_1\times (_{\xi_2} \mathcal{M}_{\gamma_2})_{X_2} & \xi_1=\gamma_1, \\
      (_{\xi_1} \mathcal{M}_{\gamma_1})_{X_1}\times \xi_2 & \xi_2=\gamma_2.
   \end{cases}
\end{equation}
Alors, le coefficient du nombre de façons de descendre dans un tel cas égale 
\begin{equation}
\label{eq4}
N_{(X_1,X_2)}(\xi,\gamma)= \begin{cases} 
     N_{X_2}(\xi_2,\gamma_2) & \xi_1=\gamma_1, \\
      N_{X_1}(\xi_1,\gamma_1) & \xi_2=\gamma_2, \\
      0 & \textnormal{autrement}.
   \end{cases}
\end{equation}
Si maintenant nous considérons l'application 
\begin{equation}
\left(\Phi: \bigoplus_{i+j=k}CM_i(M_1,f_1)\otimes CM_j(M_2,f_2)\to CM_k(M_1\times M_2, f_1+f_2)\right): \xi_1\otimes \xi_2 \mapsto (\xi_1,\xi_2).
\end{equation}
Notre objectif est de montrer que cette application est un isomorphisme de groupes abéliens. Pour ce faire, nous nous inspirerons de \cite{audin2014morse}. D'abord, considérons les lemmes, théorèmes et propositions suivants.
\begin{lemma}
L'application $\Phi$ définit un isomorphisme de complexes 
\begin{equation}
\Phi:(CM_\star(M_1,f_1)\otimes CM_*(M_2,f_2),\partial_{X_1}\otimes \mathds{1}_{M_2}+\mathds{1}_{M_1}\otimes \partial_{X_2})
\to (CM_{\star+*(M_1\times M_2, f_1+f_2)}, \partial_{X_1,X_2}).\end{equation}
\end{lemma}
\begin{prL} Par simple calcul direct, nous avons d'une part
\begin{equation}
\begin{split}
\Phi(\partial_{X_1}\otimes \mathds{1}_{M_2}+\mathds{1}_{M_1}\otimes \partial_{X_2})(\xi_1\otimes \xi_2)&=\Phi(\partial_{X_1}\xi_1\otimes \xi_2+\xi_1\otimes \partial_{X_2}\xi_2)\\&=\Phi\left(\sum_{\gamma_1\in \textnormal{Crit}_{i-1}(f_1)}N_{X_1}(\xi_1,\gamma_1)\gamma_1\otimes \xi_2+\xi_1\otimes \sum_{\gamma_2\in \textnormal{Crit}_{j-1}(f_2)}N_{X_2}(\xi_2,\gamma_2)\gamma_2\right)\\&=\sum_{\gamma_1\in \textnormal{Crit}_{i-1}(f_1)}N_{X_1}(\xi_1,\gamma_1)(\gamma_1,\xi_2)+ \sum_{\gamma_2\in \textnormal{Crit}_{j-1}(f_2)}N_{X_2}(\xi_2,\gamma_2)(\xi_1,\gamma_2).
\end{split}
\end{equation}
et d'une autre part 
\begin{equation}
\partial_{(X_1,X_2)}\Phi(\xi_1\otimes \xi_2)=\sum_{(\gamma_1,\gamma_2)\in \textnormal{Crit}_{i+j-1}(f_1+f_2)}N_{(X_1,X_2)}((\xi_1,\xi_2),(\gamma_1,\gamma_2))(\gamma_1,\gamma_2).
\end{equation}
Comme $\xi$ et $\gamma$ sont d'indices consécutifs, la somme se réécrit comme \eqref{eq4}
\begin{equation}
\partial_{(X_1,X_2)}\Phi(\xi_1\otimes \xi_2)=\sum_{\gamma_1\in \textnormal{Crit}_{i-1}(f_1)}N_{X_1}(\xi_1,\gamma_1)(\gamma_1,\gamma_2=\xi_2)+\sum_{\gamma_2\in \textnormal{Crit}_{j-1}(f_2)}N_{X_1}(\xi_1,\gamma_1)(\gamma_1=\xi_1,\gamma_2).
\end{equation}
En comparant, le résultat désiré est vérifié. Q.E.D
\end{prL}
\begin{lemma}
L'homologie d'un produit tensoriel de complexes est le produit tensoriel des homologies -- i.e. 
\begin{equation}
H_\star(C_*\otimes d_*)=H_\star(C_*)\otimes H_\star(d_*).
\end{equation}
\end{lemma}
\begin{prL}
Voir \cite[p. 556]{audin2014morse}. Q.E.D
\end{prL}
\begin{lemma}
Si $0\le \alpha <\infty$, alors 
\begin{equation}
H_\star\left(\bigoplus_\alpha C_\alpha\right)\cong \bigoplus_\alpha H_\star(C_\alpha).
\end{equation}
\end{lemma}
\begin{prL}
Voir \cite{weibel1995introduction}. Q.E.D
\end{prL}
\begin{thm} \textnormal{(Formule de Künneth)}\\

Soient $M_1$ et $M_2$ deux variétés compacte ($k$ fini). Alors, l'isomorphisme suivant existe 
\begin{equation}
HM_k(M_1\times M_2, f_1+f_2, \textnormal{Or})\to \bigoplus_{i+j=k}HM_i(M_1, f_1, \textnormal{Or})\otimes HM_j(M_2, f_2, \textnormal{Or}).
\end{equation}
\end{thm}
\begin{prT}
Par l'isomorphisme $\Phi$ (ici $\textnormal{Or}\cong \mathds{Z}$)
\begin{equation}
\begin{split}
HM_k(M_1\times M_2, f_1+f_2, \textnormal{Or})&\equiv HM(CM_k(M_1\times M_2, f_1+f_2, \textnormal{Or}))\\& \cong^{\Phi} HM\left(\bigoplus_{i+j=k} CM_i(M_1, f_1, \textnormal{Or})\otimes CM_j(M_2, f_2, \textnormal{Or}) \right)\\& \cong \bigoplus_{i+j=k} HM\left(CM_i(M_1, f_1, \textnormal{Or})\otimes CM_j(M_2, f_2, \textnormal{Or}) \right)\\& \cong  \bigoplus_{i+j=k} HM(CM_i(M_1, f_1, \textnormal{Or}))\otimes HM(CM_j(M_2, f_2, \textnormal{Or}))\\&\equiv \bigoplus_{i+j=k} HM_i(M_1, f_1, \textnormal{Or}))\otimes HM_j(M_2, f_2, \textnormal{Or}).
\end{split}
\end{equation}
Q.E.D
\end{prT}
Il est maintenant temps de définir $(f_1,X_1)$ et $(f_2,X_2)$ explicitement sur $M_1= S^2$ et sur $M_2=T^2$. Pour se faire, considérons d'abord une proposition \cite[p.9]{audin2014morse}. Elle va comme suit 
\begin{prop}
Soit $M\subset \mathds{R}^n$ une sous-variété. Pour \emph{presque} tous les points (propriété générique) $p\in \mathds{R}^n$, la fonction 
\begin{equation}
((f^p:M\to \mathds{R}): \xi \mapsto \|\xi -p\|_{\textnormal{Euc}}^2),
\end{equation}
est de Morse.
\end{prop}
\begin{prP} La différentielle de $f^p$ est donnée par 
\begin{equation}
\begin{split}
f_\star^p(\eta)&=(\frac{d}{d\xi}\|\xi-p\|^2)\cdot \eta\\&=(\frac{d}{d\xi}(\xi-p)\cdot(\xi-p))\cdot\eta\\&=2(\xi-p)\cdot \eta.
\end{split}
\end{equation}
Comme $\eta\in \textnormal{dom}(f_\star^p)\equiv \mathcal{T}_\xi M$, $\xi$ est un point critique si et seulement si $(\xi-p)$ est orthogonal à $\mathcal{T}_\xi M$. Comme $M$ est une sous-variété de $\mathds{R}^n$, par \cite[Théorème A.1.1]{audin2014morse}, il nous est possible de choisir une paramétrisation locale de $\xi\in M$ au voisinage de $p$ (si $d\le n$)
\begin{equation}
(u_1,...,u_d) \mapsto \xi(u_1,...,u_d).
\end{equation}
Dans ces coordonnées, nous avons, composantes par composante de la différentielle 
\begin{equation}
\begin{split}
\frac{\partial}{\partial_{u_i}}f^p&=\frac{\partial}{\partial_{u_i}}((\xi-p)\cdot (\xi-p))\\&=2(\xi-p)\cdot \frac{\partial}{\partial_{u_i}}\xi. 
\end{split}
\end{equation}
Alors, 
\begin{equation}
\begin{split}
\frac{\partial^2}{\partial_{u_i}\partial_{u_j}}f^p&=\frac{\partial}{\partial_{u_i}}(2(\xi-p)\cdot \frac{\partial}{\partial_{u_j}}\xi)\\&=2\left( \frac{\partial}{\partial_{u_i}\partial_{u_i}}\xi\cdot \frac{\partial}{\partial_{u_i}\partial_{u_j}}\xi+(\xi-p)\cdot \frac{\partial^2}{\partial_{u_i}\partial_{u_j}}\xi\right). 
\end{split}
\end{equation}
Rappelons que si $A\in \textnormal{MAT}_{n\times n}(V,\mathds{K})$, $A$ est inversible (non dégénérée) si et seulement si $\det(A)\neq 0$ si et seulement si $\textnormal{Rang}(A)=n$ si et seulement si $\ker(A)=0$. Ainsi, le point $\xi$ est un point critique non dégénéré si et seulement si le vecteur $(\xi-p)$ est orthogonal à $\mathcal{T}_\xi M$ et que le rang de $\left[\frac{\partial^2}{\partial_{u_i}\partial_{u_j}}f^p\right]$ est égal à $d$. Il reste à montrer que, sous ces conditions, $f^p$ est génériquement une fonction de Morse. Il est suffisant de montrer que les $p$ qui ne satisfont pas cette condition de générécité sur la non-dégénérescence sont des points critiques d'une application lisse dans le but d'utiliser le Théorème de Sard.\\

Pour ce faire, nous considèrerons le fibré normal de $M$ dans $\mathds{R}^n$ -- i.e. 
\begin{equation}
N:=\{(\xi,\nu)\in M\times \mathds{R}^n \mid \nu \perp \mathcal{T}_\xi M\}\subset M\times \mathds{R}^n,
\end{equation}
muni de l'application 
\begin{equation}
(E:N\to \mathds{R}^n):(\xi,\nu)\mapsto \xi+\nu.
\end{equation}
Le résultat désiré suit du prochain lemme. Q.E.D
\end{prP}
En effet, nous avons montré plus haut que $f^p$ était de Morse à $\xi$ (point critique non dégénéré) si et seulement si la matrice  $\left[\frac{\partial^2}{\partial_{u_i}\partial_{u_j}}f^p\right]$ est inversible à ce point. Ainsi, $f^p$ n'est pas de Morse à $\xi$ si et seulement si la matrice  $\left[\frac{\partial^2}{\partial_{u_i}\partial_{u_j}}f^p\right]$ est non-inversible à ce point. Le lemme suivant nous garantit que si cette matrice n'est pas inversible, alors, $p$ est un point critique de $E$ (application lisse). Par le Théorème de Sard, $\textnormal{Crit}(E)$ constitue un ensemble fini de points. De plus, comme $p:=\xi+\nu$, il existe clairement une bijection du point $\xi$ à $p$. Ainsi, il doit exister un nombre fini de $\xi$ qui ne sont pas des points critiques non dégénérés de $f^p$ -- i.e. $f^p$ est de Morse presque partout. 
\begin{lemma}
Le fibré normal $N$ est une sous-variété de $M\times \mathds{R}^n$. Le point $p:=\xi+\nu\in \mathds{R}^n$ est une valeur critique de $E$ si et seulement si, la matrice avec les éléments donnés par 
\begin{equation}
\begin{split}
\frac{\partial^2}{\partial_{u_i}\partial_{u_j}}f^p&=2\left( \frac{\partial}{\partial_{u_i}\partial_{u_i}}\xi\cdot \frac{\partial}{\partial_{u_i}\partial_{u_j}}\xi-\nu\cdot \frac{\partial^2}{\partial_{u_i}\partial_{u_j}}\xi\right),
\end{split}
\end{equation}
n'est pas inversible.
\end{lemma}
\begin{prL}
Voir \cite[p.10]{audin2014morse}. Q.E.D
\end{prL}
De cette dernière proposition, soit $p\not\in S^2,T^2$ et leur intérieur. Alors, il est facile d'obtenir la fonction de Morse suivante définie sur $S^2\times T^2$
\begin{equation}
(F_p:=f_1^p+f_2^p: S^2\times T^2\to \mathds{R}): \|\xi-p\|^2+\|\eta-p\|^2.
\end{equation}
Ainsi, prenons $X_1:=\Grad(f_1^p)$ et $X_2:=\Grad(f_2^p)$ au sens de \cite{audin2014morse}, afin que $(X_1,X_2)$ soit de Smale sur $S^2\times T^2$.\\

Nous devons maintenant construire le complexe de Morse pour $S^2$ et $T^2$ et ensuite utiliser la formule de Künneth pour calculer les groupes d'homologie du produit. \\

Le complexe pour $S^2$ s'écrit, comme cet espace n'est muni que d'un maximum $\alpha$ d'indice 2 et un minimum $\beta$ d'indice $0$,
\begin{equation}
...\xrightarrow{\partial_4}\bigoplus_{x_3\in \textnormal{Crit}_3(f_1^p)}C_{x_3}x_3\xrightarrow{\partial_3}\bigoplus_{x_2\in \textnormal{Crit}_2(f_1^p)}C_{x_2}x_2\xrightarrow{\partial_2} \bigoplus_{x_1\in \textnormal{Crit}_1(f_1^p)}C_{x_1}x_1 \xrightarrow{\partial_1} \bigoplus_{x_0\in \textnormal{Crit}_0(f_1^p)}C_{x_0}x_0  \xrightarrow{\partial_0} 0,
\end{equation}
ou 
\begin{equation}
...\xrightarrow{\partial_4}0\xrightarrow{\partial_3}C_{\alpha}\alpha\xrightarrow{\partial_2} 0\xrightarrow{\partial_1} C_{\beta}\beta  \xrightarrow{\partial_0} 0,
\end{equation}
ou, par isomorphisme,
\begin{equation}
...\xrightarrow{\partial_4}0\xrightarrow{\partial_3}\mathds{Z}\xrightarrow{\partial_2} 0\xrightarrow{\partial_1} \mathds{Z} \xrightarrow{\partial_0} 0.
\end{equation}
Il est aussi à noter que l'application de bord $\partial_2$ est telle que 
\begin{equation}
\partial_2\alpha=\sum_{y\in \textnormal{Crit}_1(f_1^p)}N_{X_1}(\alpha,y)y=0,
\end{equation}
puisque la somme est vide. Il en va de même pour $\partial_0\beta$. Les groupes d'homologies sont donc, pour $S^2$, 
\begin{equation}
HM_\star(S^2,f_1^p,\mathds{Z})= \begin{cases} 
      \mathds{Z} & \star\in\{0,2\}, \\
       0 & \textnormal{autrement}.
   \end{cases}
\end{equation}
De façon similaire, il nous est possible de construire le complexe de Morse pour $T^2$, gardant en tête que cet espace admet avec notre choix de $f_2^p$, un maximum $a$ d'indice $2$, deux points de scelle $c_1$ et $c_2$ d'indice 1 et finalement un minimum $b$ d'indice $0$. Ainsi, 
\begin{equation}
...\xrightarrow{\partial_4}0\xrightarrow{\partial_3}\mathds{Z}\xrightarrow{\partial_2} \mathds{Z}\oplus \mathds{Z}\xrightarrow{\partial_1} \mathds{Z} \xrightarrow{\partial_0} 0.
\end{equation}
Ici, les applications de bord sont 
\begin{equation}
\partial_2 a=\sum_{z\in \textnormal{Crit}_1(f_2^p)}N_{X_2}(a,z)z=\sum_{z\in \textnormal{Crit}_1(f_2^p)}\sum_{u\in _a \mathcal{M}_z}n_u z=((+1)+(-1))c_1+((+1)+(-1))c_2=0.
\end{equation}
\begin{equation}
\partial_1 c_1=\sum_{l\in \textnormal{Crit}_0(f_2^p)}N_{X_2}(c_1,l)l=\sum_{l\in \textnormal{Crit}_0(f_2^p)}\sum_{v\in _{c_1} \mathcal{M}_l}n_v z=((+1)+(-1))b=0=\partial_1 c_2.
\end{equation}
Ainsi, les groupes d'homologie sont (sans surprise), pour $T^2$,
\begin{equation}
HM_\star(T^2,f_2^p,\mathds{Z})= \begin{cases} 
      \mathds{Z} & \star\in\{0,2\}, \\
      \mathds{Z}\oplus \mathds{Z} & \star\in\{1\}, \\
       0 & \textnormal{autrement}.
   \end{cases}
\end{equation}
À l'aide de la formule de Künneth, calculons les groupes d’homologie. Pour $k\ge 5$, c'est groupes sont triviaux par les groupes d'homologie trouvés pour $S^2$ et $T^2$. Ainsi, pour $k=4$, 
\begin{equation}
\begin{split}
HM_4(S^2\times T^2, f_1^p+f_2^p,\mathds{Z})&\cong \bigoplus_{i+j=4}HM_i(S^2, f_1^p,\mathds{Z})\otimes HM_j(T^2, f_2^p,\mathds{Z})\\&=HM_0(S^2, f_1^p,\mathds{Z})\otimes HM_4(T^2, f_2^p,\mathds{Z})\oplus HM_1(S^2, f_1^p,\mathds{Z})\otimes HM_3(T^2, f_2^p,\mathds{Z}) \\& \oplus HM_2(S^2, f_1^p,\mathds{Z})\otimes HM_2(T^2, f_2^p,\mathds{Z}) \oplus HM_3(S^2, f_1^p,\mathds{Z})\otimes HM_1(T^2, f_2^p,\mathds{Z}) \\ & \oplus HM_4(S^2, f_1^p,\mathds{Z})\otimes HM_0(T^2, f_2^p,\mathds{Z})\\&= \mathds{Z}\otimes \mathds{Z}.
\end{split}
\end{equation}
Pour $k=3$ 
\begin{equation}
\begin{split}
HM_3(S^2\times T^2, f_1^p+f_2^p,\mathds{Z})&\cong \bigoplus_{i+j=3}HM_i(S^2, f_1^p,\mathds{Z})\otimes HM_j(T^2, f_2^p,\mathds{Z})\\&=HM_0(S^2, f_1^p,\mathds{Z})\otimes HM_3(T^2, f_2^p,\mathds{Z})\oplus HM_1(S^2, f_1^p,\mathds{Z})\otimes HM_2(T^2, f_2^p,\mathds{Z}) \\& \oplus HM_2(S^2, f_1^p,\mathds{Z})\otimes HM_1(T^2, f_2^p,\mathds{Z}) \oplus HM_3(S^2, f_1^p,\mathds{Z})\otimes HM_0(T^2, f_2^p,\mathds{Z})\\&= \mathds{Z}\otimes( \mathds{Z}\oplus \mathds{Z}).
\end{split}
\end{equation}
Pour $k=2$
\begin{equation}
\begin{split}
HM_2(S^2\times T^2, f_1^p+f_2^p,\mathds{Z})&\cong \bigoplus_{i+j=2}HM_i(S^2, f_1^p,\mathds{Z})\otimes HM_j(T^2, f_2^p,\mathds{Z})\\&=HM_0(S^2, f_1^p,\mathds{Z})\otimes HM_2(T^2, f_2^p,\mathds{Z})\oplus HM_1(S^2, f_1^p,\mathds{Z})\otimes HM_1(T^2, f_2^p,\mathds{Z}) \\& \oplus HM_2(S^2, f_1^p,\mathds{Z})\otimes HM_0(T^2, f_2^p,\mathds{Z})\\&= (\mathds{Z}\otimes \mathds{Z})\oplus (\mathds{Z}\otimes \mathds{Z}).
\end{split}
\end{equation}
Pour $k=1$
\begin{equation}
\begin{split}
HM_1(S^2\times T^2, f_1^p+f_2^p,\mathds{Z})&\cong \bigoplus_{i+j=1}HM_i(S^2, f_1^p,\mathds{Z})\otimes HM_j(T^2, f_2^p,\mathds{Z})\\&=HM_0(S^2, f_1^p,\mathds{Z})\otimes HM_1(T^2, f_2^p,\mathds{Z})\oplus HM_1(S^2, f_1^p,\mathds{Z})\otimes HM_0(T^2, f_2^p,\mathds{Z})\\&= \mathds{Z}\otimes( \mathds{Z}\oplus \mathds{Z}).
\end{split}
\end{equation}
Pour $k=0$
\begin{equation}
\begin{split}
HM_1(S^2\times T^2, f_1^p+f_2^p,\mathds{Z})&\cong \bigoplus_{i+j=0}HM_i(S^2, f_1^p,\mathds{Z})\otimes HM_j(T^2, f_2^p,\mathds{Z})\\&=HM_0(S^2, f_1^p,\mathds{Z})\otimes HM_0(T^2, f_2^p,\mathds{Z})\\&= \mathds{Z}\otimes \mathds{Z}.
\end{split}
\end{equation}
En récapitulatif, les groupes d'homologie sont 
\begin{equation}
HM_\star(S^2\times T^2,f_1^p+f_2^p,\mathds{Z})\cong \begin{cases} 
      \mathds{Z}\otimes \mathds{Z} & \star\in\{0,4\}, \\
     \mathds{Z}\otimes (\mathds{Z}\oplus \mathds{Z}) & \star\in\{1,3\}, \\
      (\mathds{Z}\otimes\mathds{Z})\oplus (\mathds{Z}\otimes\mathds{Z}) & \star\in\{2\}, \\
       0 & \textnormal{autrement}.
   \end{cases}
\end{equation}
Comme le produit tensoriel est distributif, nous pouvons réécrire 
\begin{equation}
HM_\star(S^2\times T^2,f_1^p+f_2^p,\mathds{Z})\cong \begin{cases} 
      \mathds{Z}\otimes \mathds{Z} & \star\in\{0,4\}, \\
     (\mathds{Z}\otimes\mathds{Z})\oplus (\mathds{Z}\otimes\mathds{Z}) & \star\in\{1,2,3\}, \\
       0 & \textnormal{autrement}.
   \end{cases}
\end{equation}
De plus, notons que pour un indice $J$ dénombrant des groupes abéliens $A_j, \ j\in J$, nous avons sur les dimensions des groupes 
\begin{equation}
\dim\left(\bigoplus_{j\in J}A_j\right)=\sum_{j\in J}\dim(A_j),
\end{equation}
et 
\begin{equation}
\dim(A_1\otimes A_2)=\dim(A_1)\dim(A_2).
\end{equation}
Comme $\mathds{Z}$ est un groupe cyclique, il s'en suit que 
\begin{equation}
\dim(HM_\star(S^2\times T^2,f_1^p+f_2^p,\mathds{Z}))= \begin{cases} 
      1 & \star\in\{0,4\}, \\
     2& \star\in\{1,2,3\}, \\
       0 & \textnormal{autrement}.
   \end{cases}
\end{equation}
Comme $(\xi_1,\xi_2)$ est un point critique de $f_1^p+f_2^p$ si et seulement si $\xi_1$ et $\xi_2$ sont des points critiques de $f_1^p$ et $f_2^p$ respectivement, on voit directement que l'égalité à l'inégalité de Morse nous donne le nombre de points critiques sur $S^2\times T^2$ pour chaque indice : 1 point critique d'indice 0, $z$, et d'indice 4, $Q$, et 2 points critiques d'indice 1, $u_1$ et $u_2$, d'indice 2, $d_1$ et $d_2$, et d'indice 3, $T_1$ et $T_2$. Le polynôme de Poincaré 
\begin{equation}
P_{S^2\times T^2}(t)= 1+2(t+t^2+t^3)+t^4,
\end{equation}
évalué à $t=-1$ nous donne la caractéristique d'Euler de $S^2\times T^2$, soit $\chi(S^2\times T^2)=0$. Ainsi, $S^2\times T^2$ est de genre 1. 
\subsection{L'éclatement des points critiques de $S^2\times T^2$ }
Soit $(f,X)$ un couple de Morse-Smale sur une variété compacte $M$. De \cite{audin2014morse}, nous avons que 
\begin{lemma} Étant donné une fonction de Morse $h$ sur une variété compacte $M$, il existe une seconde fonction de Morse $\ell$ telle que $\textnormal{Crit}_k(\ell)=\textnormal{Crit}_k(h) \ \ \forall \ k$ et telle que $\ell(p)=|p|$.
\end{lemma}
Notons qu'une telle fonction de Morse est nommée \emph{ordonnée}.
\begin{lemma}
Sur toute variété compacte et connexe $M$, il existe une fonction de Morse qui ne possède qu'un seul minimum. 
\end{lemma}
Ainsi, supposons que la fonction $f$ est ordonnée et possède un unique minimum, noté $z$.\\

Soit maintenant $p\in \textnormal{Crit}(f)$. Rappelons que dans \cite{audin2014morse}, on établie que $\overline{_p \mathcal{M}_z}$ est la compactification de l'espace des modules de flots reliant $p$ à $z$ par l'ajout (ou non) des trajectoires brisées.\\

Ainsi, pour $v\in \overline{_p \mathcal{M}_z}$, on peut écrire $v=(v_1,...,v_j)\in\prod_{i\in \{1,...,j\}, p_1=p, p_j=z} \overline{_{p_i} \mathcal{M}_{p_{i+1}}}$, tel que $p_i\in \textnormal{Crit}(f)$.\\

Il sera utile de paramétrer ces trajectoires de la façon la plus naturelle. Pour $v\in _p\mathcal{M}_z$, choisissons la paramétrisation
\begin{equation}
((\Lambda_v:[0,f(p)]\to M): \Lambda_v(\tau)=\zeta \Leftrightarrow f(\zeta)=f(p)-\tau).
\end{equation}
Moralement, cette paramétrisation nous informe du \emph{niveau} $\zeta$ sur lequel on se trouve après un \emph{temps} $\tau$ compris dans l'intervalle du domaine. Cette paramétrisation s'étend donc bel et bien à toutes les trajectoires de modules de flots (brisées ou non) dans $\partial \overline{_p \mathcal{M}_z}$, ayant comme limites évidentes 
\begin{equation}
\Lambda_v(\tau)\xrightarrow{\tau \to 0} p \ \ \ \textnormal{et} \ \ \ \Lambda_v(\tau)\xrightarrow{\tau \to f(p)} 0.
\end{equation}
\begin{definition} \textnormal{(Éclatement d'une variété instable)} L'éclatement d'une variété instable est par définition 
\begin{equation}
\mathcal{E}_M(p,\tau): (\overline{_p \mathcal{M}_z}\times [0,f(p)])/\sim_\tau, 
\end{equation}
où la relation d'équivalence $\sim_\tau$ est telle que $((v_1,...,v_j),\tau) \sim_\tau ((v_1',...,v_j'),\tau)$ si, $\forall \ i : f(p_{i-1})>\tau$, on a $v_i=v_i'\in _{p_{i-1}}\mathcal{M}_{p_i}$.
\end{definition}
Autrement dit, on identifie les points $(v,\tau)$ et $(v',\tau)$ si les trajectoires de $v$ et de $v'$ coïncident partout au-dessus du niveau $\tau$.\\

Avant de tenter de visualiser l'éclatement des points critiques de $S^2\times T^2$, voici quelques résultats intéressants sur l'espace $\mathcal{E}_M(p,t)$. 
\begin{lemma}
$\mathcal{E}_M(p,t)$ est Hausdorff.
\end{lemma}
\begin{prL}
 Rappelons qu'un espace est Hausdorff s’il existe pour tous les points distincts un voisinage pour chaque ne s'intersectant pas. Autrement dit, que n'importe quelle suite qui converge dans cet espace est unique. Ainsi, soient deux suites $(v_n,\tau_n)$ et $(v_n',\tau_n')$ dans $\overline{_p \mathcal{M}_z}\times [0,f(p)]$ telles que : (i) $v_n \xrightarrow{n\to \infty}v$ et $\tau_n \xrightarrow{n\to \infty}\tau$, (ii) $v'_n \xrightarrow{n\to \infty}v'$ et $\tau_n' \xrightarrow{n\to \infty}\tau'$ et (iii)  $(v_n,\tau_n)|_\sim=(v_n',\tau_n')|_\sim \  \forall \ n$.\\
 
 On doit alors montrer que ces deux suites de trajectoires, donnée la relation d'équivalence, sont égales dans la limite où $n\to \infty$.\\
 
D'abord, comme $(v_n,\tau_n)\sim(v_n',\tau_n') \  \forall \ n$, par définition de $\sim_\tau$, on doit avoir que $\tau_n=\tau_n'$ pour tout $n$. Ensuite, comme $(v_n,\tau)\sim(v_n',\tau) \  \forall \ n$, nous avons que passé le temps limite $\tau$ (ou, au-dessus du niveau $\tau$), toutes les trajectoires $v_n$ et $v_n'$ coïncident. Q.E.D
\end{prL}
La paramétrisation $\Lambda_v$ nous permet de définir l'application continue 
\begin{equation}
((\Lambda^*:\overline{_p \mathcal{M}_z}\times [0,f(p)]\to M):(v,\tau)\mapsto \Lambda_v(f(p)-\tau)=\iota \Leftrightarrow f(\iota)=\tau). 
\end{equation}
Cette application permet d'identifier un point $\iota$ sur la variété à l'intersection de $v$ et du niveau $f^{-1}(\tau)$ -- i.e. $\iota\in W^u(p)$ ou sa frontière. Comme le couple $(v,\tau)$ est unique, cette identification est injective.\\

L'application $\Lambda^*$ peut être redéfinie par 
\begin{equation}
\Lambda^*\equiv E\circ \sim_\tau,
\end{equation}
où 
\begin{equation}
((E:\mathcal{E}_M(p,\tau)\to M):(v,\tau)|_\sim\mapsto \Lambda^*(v,\tau)=\Lambda_v(f(p)-\tau)).
\end{equation}
Le lemme suivant sur l'espace engendré par l'image de $E$ justifie pourquoi on appelle $\mathcal{E}_M(p,\tau)$ l'éclatement de $W^u(p)$.
\begin{lemma}
$E(\mathcal{E}_M(p,\tau))=\overline{W^u(p)}$. 
\end{lemma}
\begin{prL}
Montrons d'abord que $E(\mathcal{E}_M(p,\tau))\supseteq \overline{W^u(p)}$. Prenons $\xi \in W^u(p)$. Alors, $\xi$ appartient au premier segment d'une trajectoire brisée (ou non) $v=(v_1,...,v_j)$ entre $p$ et $z$ (ici, $j\le 1$ est situé au-dessus d'un certain niveau $f^{-1}(\tau)$). Ainsi, par définition de $\Lambda^*$ et donc de $E$, il existe un couple $(v,\tau)$ tel que $E((v,\tau)|_\sim)=\Lambda^*((v,\tau))=\xi$. Ainsi, pour tout $\xi \in W^u(p)$, $E(\mathcal{E}_M(p,\tau))\supseteq W^u(p)$. Pour montrer totalement l'inclusion désirée, nous devons prouver la même chose pour $\xi \in \partial W^u(p)$. Soit alors $(\xi_n)_{n=1}^{\infty}\subseteq \overline{W^u(p)}$ une séquence telle qu'elle converge à $\xi$ (une telle séquence existe car $\overline{W^u(p)}=W^u(p)\cup \partial W^u(p)$ est fermé et donc contient tous ses points d'accumulation et que $\xi$ en est trivialement un). Par l'argument développé plus haut, pour tout $n$, on trouve un $(v_n,\tau_n): \Lambda^*((v_n,\tau_n))=\xi_n$. Comme $((v_n,\tau_n))_{n=1}^\infty$ est une suite définie sur un compact, par Bolzano-Weierstra\ss \ généralisé, il existe une sous-suite $((v_{n_k},\tau_{n_k}))_{k_1}^{\infty}$ qui converge dans ce compact. Ainsi, il existe $(v,\tau)\in \overline{_p \mathcal{M}_z}$ tel que, par continuité,  
\begin{equation}
\Lambda^*((v_{n_k},\tau_{n_k}))=\xi_{n_k}\xrightarrow{k\to \infty} \Lambda^*((v,\tau)).
\end{equation}
Comme l'éclatement est Hausdorff, la limite des suites dans cet espace est unique  et comme $\xi_n\xrightarrow{n\to \infty} \xi$, alors
\begin{equation}
E((v_{n_k},\tau_{n_k})|_\sim)\xrightarrow{k\to \infty} E((v,\tau)|_\sim)\equiv (E\circ \sim)((v,\tau))\equiv \Lambda^*((v,\tau))=\xi.
\end{equation}
La première inclusion est ainsi montrée.\\

Il faut maintenant montrer que $E(\mathcal{E}_M(p,\tau))\subseteq \overline{W^u(p)}$. Pour se faire, prenons $\gamma\in E(\mathcal{E}_M(p,\tau))$. Alors, $\gamma$ est tel que 
\begin{equation}
(\gamma =(v,\tau)|_\sim \mid (v,\tau)\in \overline{_p \mathcal{M}_z}).
\end{equation}
Il faut montrer que $E(\gamma)\in \overline{W^u(p)}$. Si $(v,\tau)\in \overline{_p \mathcal{M}_z}\times (0, f(p)]$, la preuve est triviale car ce point n'est pas un point limite -- i.e. on s'y rend dans un temps non asymptotique. C'est-à-dire qu'il existe un $\gamma$ bien défini tel que $E(\gamma)\in W^u(p)\subseteq \overline{W^u(p)}$. Si $(v,\tau)\not\in \overline{_p \mathcal{M}_z}\times (0, f(p)]$, alors, choisissons une séquence $(v_n,\tau_n)_{n=1}^\infty\subset \overline{_p \mathcal{M}_z}\times (0, f(p)] : (v_n,t_n)\xrightarrow{n\to \infty}(v,t)$\footnote{Clairement, $(v,t)$ est un point d'accumulation de $\overline{_p \mathcal{M}_z}\times (0, f(p)]$ données $v_n\xrightarrow{n\to \infty}v$ et $\tau_n\xrightarrow{n\to \infty}\tau$.}. Par la continuité de $\Lambda^*$, alors $\Lambda^*((v_n,\tau_n))\xrightarrow{n\to \infty}\Lambda^*((v,\tau))\equiv (E\circ \sim)((v,\tau))\equiv E(\gamma)$. Alors, $E(\gamma)= \Lambda^*((v,\tau)) \in \overline{W^u(p)}$. Q.E.D
\end{prL}
Il est intuitif que l'éclatement soit un espace contractible -- i.e.  l'application identité y est nulle-homotopique à une application constante.
\begin{lemma}
L'espace $\mathcal{E}_M(p,\tau)$ est contractible. 
\end{lemma}
\begin{prL}
Par définition, les points de l'espace $\overline{_p \mathcal{M}_z}\times \{f(p)\}$ sont envoyés sur le même point dans l'espace de l'éclaté -- par le quotient qui le défini. Notons ce point $\varpi\in \mathcal{E}_M(p,\tau)$. Comme dit plus tôt, nous désirons montrer que $\textnormal{Id}(\mathcal{E}_M(p,\tau))$ est homotopique à une application constante envoyant $\mathcal{E}_M(p,\tau)$ sur $\varpi$ -- i.e. application dont l'image est $\varpi$.\\

Pour se faire, définissons 
\begin{equation}
((\Lambda_v':[0,f(p)]\to \mathcal{E}_M(p,\tau)):\tau\mapsto (v,f(p)-\tau)|_\sim).
\end{equation}
Alors, $\Lambda_v'(0)=(v, f(p))|_\sim=\varpi$ et $\Lambda_v\equiv E\circ \Lambda_v'$.\\

Munis de cette notation, on note qu'il y a pour tout $\zeta\in \mathcal{E}_M(p,\tau)$ une trajectoire $\bar{\zeta}:[0,\tau^*]\to \mathcal{E}_M(p,\tau)$ telle que $\bar{\zeta}(0)=\varpi$ et $\bar{\zeta}(\tau^*)=\zeta$ qui, de plus, coïncide avec une trajectoire $v\in \overline{_p \mathcal{M}_z}$ -- i.e.
\begin{equation}
\bar{\zeta}(\tau)=\Lambda_v'(\tau) \ \forall \ \tau\in [0,\tau^*].
\end{equation}
Pour une trajectoire $\bar{\zeta}$, $\tau^*$ se doit d'être unique. On note que l'association constante, $\zeta \to \bar{\zeta}$,
\begin{equation}
((\beta:\mathcal{E}_M(p,\tau)\to \underbrace{\{\bar{\zeta}:[0,\tau^*]\to \mathcal{E}_M(p,\tau) \mid \tau^*\ge 0, \bar{\zeta}(0)=\varpi\}}_{\textnormal{Ensemble des trajectoires $\varpi\to \zeta$}}):\zeta \mapsto \bar{\zeta}).
\end{equation}
Alors, il nous est possible de définir l'homotopie $\mathcal{H}_\alpha$ 
\begin{equation}
((\mathcal{H}_\alpha: \mathcal{E}_M(p,\tau)\times [0,1]\to \mathcal{E}_M(p,\tau)):\zeta \mapsto (\beta(\zeta))((1-\alpha)\tau^*)\equiv \bar{\zeta}((1-\alpha)\tau^*)), 
\end{equation}
telle que 
\begin{align}
\mathcal{H}_0(\zeta)=\bar{\zeta}(\tau^*)=\textnormal{Id}(\zeta),\\
 \mathcal{H}_1(\zeta)=\bar{\zeta}(0)=\varpi.
\end{align}
Il existe donc une homotopie entre l'identité et une application constante. Ainsi, l'identité est nulle-homotopique. Q.E.D 
\end{prL}
\begin{lemma}
$\mathcal{E}_M(p,\tau)$ est homéomorphique à un disque fermé de dimension égale à l'indice de $p$. De plus, 
\begin{equation}
\partial \mathcal{E}_M(p,\tau)=\bigcup_{\textnormal{Ind}(q)<\textnormal{Ind}(p)} \overline{_p\mathcal{M}_q}\times \mathcal{E}_M(q,\tau).
\end{equation}
\end{lemma}
\begin{prL}
Voir \cite[Section 2.4.6]{barraud2007lagrangian}. Q.E.D
\end{prL}
Effectuons maintenant l'éclaté du point $T_1$. Comme l'indice décroît le long des trajectoires d'un pseudo-gradient satisfaisant la condition de Smale, lorsque nous utilisons la définition des trajectoires brisées
\begin{equation}
\overline{ _\xi \mathcal{M}_\gamma}=\bigsqcup_{\substack{\sigma\in \textnormal{Crit}_h(f), \\ \textnormal{Ind}(\xi)>h>\textnormal{Ind}(\gamma)}}  {_\xi\mathcal{M}_\sigma} \times  {_\sigma\mathcal{M}_\gamma},
\end{equation}
il ne faut impérativement pas considérer les trajectoires entre points du même indice.\\

$\overline{ _{T_1} \mathcal{M}_z}$ par définition, exprimé en terme de fermeture de cellules $e^k$
\begin{equation}
\begin{split}
\overline{ _{T_1} \mathcal{M}_z}&=\bigsqcup_{\substack{\sigma\in \textnormal{Crit}_h(f), \\ 3>h>0}}  {_\xi\mathcal{M}_\sigma} \times  {_\sigma\mathcal{M}_\gamma}\\&= \{{_{T_1}\mathcal{M}_{d_1}}\times {_{d_1}\mathcal{M}_{u_1}}\times {_{u_1}\mathcal{M}_{z}}\}\cup \{{_{T_1}\mathcal{M}_{d_1}}\times {_{d_1}\mathcal{M}_{u_2}}\times {_{u_2}\mathcal{M}_{z}}\}\\&\cup \{{_{T_1}\mathcal{M}_{d_2}}\times {_{d_2}\mathcal{M}_{u_1}}\times {_{u_1}\mathcal{M}_{z}}\}\cup \{{_{T_1}\mathcal{M}_{d_2}}\times {_{d_2}\mathcal{M}_{u_2}}\times {_{u_2}\mathcal{M}_{z}}\}\\& \cup \{{_{T_1}\mathcal{M}_{u_1}}\times {_{u_1}\mathcal{M}_{z}}\}\cup \{{_{T_1}\mathcal{M}_{u_2}}\times {_{u_2}\mathcal{M}_{z}}\}\\&\cup \{{_{T_1}\mathcal{M}_{d_1}}\times {_{d_1}\mathcal{M}_{z}}\} \cup \{{_{T_1}\mathcal{M}_{z}}\}\\&= \bigcup_{i=1}^4(\{\bar{e}^0, \bar{e}^0\}_i\cup \{\bar{e}^1\}_i)\cup \{\bar{e}^2\}.
\end{split}
\end{equation}
Alors, l'éclaté de $T_1$ est 
\begin{equation}
\mathcal{E}_{S^2\times T^2}(T_1,\tau)= \overline{ _{T_1} \mathcal{M}_z}\times [0,f(T_1)=3]=\bigcup_{i=1}^4(\{\bar{e}^1, \bar{e}^1\}_i\cup \{\bar{e}^2\}_i)\cup \{\bar{e}^3\}.
\end{equation}



\section{Invitation à une généralisation}
Pour conclure, nous mentionnons que dans le futur, il serait intéressant de transporter la discussion faite ici, mais dans le cas où la dimension est infinie. Elle est définie sur un espace de
dimension infini, comme par exemple celui des lacets d'une variété symplectique, avec la fonctionnelle d'action
comme fonction de Morse. Bien que l'indice et le coindice en tout point critique soient infinis, on peut
donner un sens à la notion d'indice relatif. Cela donne lieu à la fameuse homologie de Floer \cite{floer1988instanton, floer1988morse, floer1988unregularized, floer1989symplectic}.



\newpage
\appendix
\section{Courte révision de l'homologie de Morse sur $\mathds{Z}$}
Le complexe de Morse associé à la fonction de Morse $f$, avec coefficients entiers et indices de Morse gradées, forme une chaîne de groupes abéliens libres générés par les points critiques dans $\textnormal{Crit}_k(f)$ -- i.e.
\begin{equation}
CM_k(M,f):=\bigoplus_{\xi\in \textnormal{Crit}_k(f)}C_\xi \xi, \ \ k,C_\xi\in \mathds{Z}.
\end{equation}
Une somme sur $\emptyset$ est par définition triviale.\\

Suivant le même algorithme, nous choisissons d'abord une orientation pour chaque variété instable. L'ensemble de ces choix est dénoté: $\textnormal{Or}$. De plus, dans ce travail, nous  assumons que $\textnormal{Ind}(\xi)-\textnormal{Ind}(\gamma)=1$
et que $u\in _\xi \mathcal{M}_\gamma$. Dans ce cas, l'espace des modules de flots $_\xi \mathcal{M}_\gamma$ est une variété orientée compacte de dimension $0$ et donc correspond à un nombre fini de points, chacun muni d'un signe $\pm 1$. Particulièrement, cet espace peut 
être associé à l'intersection des 0-sections d'un fibré vectoriel de Banach \cite[p.125]{fell1988representations}. L'orbite $_\xi M_\gamma^u$ est la composante connectée de $_\xi M_\gamma$ passant par $u$, et est donc munie d'une orientation (par le truchement de la première proposition) $\mathcal{O}(_\xi M_\gamma^u)$. L'orientation induite par le flot $\mathcal{O}^{\dot{\phi}_t}(u)$ et le signe caractéristique $\eta_u:=\eta_u(\textnormal{Or})$ est définie par la relation 
\begin{equation}
\mathcal{O}(_\xi M_\gamma^u)=\eta_u \mathcal{O}^{\dot{\phi}_t}(u).
\end{equation}
Il nous est alors possible de définir l'opérateur de bord 
\begin{equation}
\left(\left(\partial_k:=\partial_k(M,f,\textnormal{Or}):CM_k(M,f)\to CM_{k-1}(M,f)\right):\xi\mapsto \sum_{\gamma\in \textnormal{Crit}_{k-1}f}N(\xi,\gamma)\gamma\right),
\end{equation}
où 
\begin{equation}
N(\xi,\gamma):=\sum_{u\in _\xi \mathcal{M}_\gamma}\eta_u, \ \  \eta_u\in \{\pm 1\} \ \textnormal{défini plus haut}.
\end{equation}
Comme cet opérateur est integrable \cite{latour1994existence}, $\partial_k^2\equiv 0$, nous nous avons que le $k$-ième groupe d'homologie correspond alors au co-noyau 
\begin{equation}
HM_k(M,f,\textnormal{Or})=\frac{\ker(\partial_k)}{\textnormal{Im}(\partial_{k+1})}.
\end{equation}
\section{Discussion sur l'existence de la limite}\label{appC}
Rappelons que $_\xi \mathcal{M}_\gamma$ représente l'ensemble des orbites du flot du pseudo gradient $\phi_\eta$ de $\xi$ à $\gamma$. Si nous avons $u\in _\xi \mathcal{M}_\gamma$ comme point de départ, pour $\eta\in [\eta_0,\infty[$, $\phi_\eta(u)\in _\xi \mathcal{M}_\gamma$ impliquant $\dot{\phi}_\eta(u)\in T _\xi \mathcal{M}_\gamma$. De \cite[Lemme 8.5]{cohen1991topics}, il s'en suit que
\begin{lemma} Soit $H(t)$ un opérateur égalant la Hessienne de $f$ dans la limite $t\to \pm \infty$. L'équation 
\begin{equation}
\frac{\textnormal{d}}{\textnormal{d}t} \dot{\phi}_t(u)+H(t)\dot{\phi}_t(u)=0,
\end{equation}
admet une solution pour $|t|>T$ telle que, $\dot{\phi}_t(u)\xrightarrow{t\to \pm \infty}0$ -- i.e. on atteint les points critiques dans cette limite et les tangentes au flot y sont nulles.  
\end{lemma}
Si nous définissons l'opérateur différentiel $\mathfrak{d}:=\frac{\textnormal{d}}{\textnormal{d}t}+H(t)$, alors par \cite[Lemme B.5]{schwarz1993morse}, comme 
\begin{equation}
\dot{\phi}_t(u)\in \ker(\mathfrak{d}),
\end{equation}
alors, gardant en tête que dans la limite $H(t)$ est non dégénérée comme $f$ est une fonction de Morse,
\begin{equation}
\left(\left(\lim_{t\to \pm \infty}\frac{\dot{\phi}_t(u)}{\|\dot{\phi}_t(u)\|_{\textnormal{Euc.}}}:= \mathfrak{E}^{\pm} \ \  \exists\right):H(\pm \infty)\mathfrak{E}^{\pm}=E^\pm \mathfrak{E}^{\pm}, E^+>0 \ \wedge \ E^- <0 \right),
\end{equation}
où il est naturel de définir $\mathfrak{E}^+:=\star^{\dot{\phi}_{t\to + \infty}}(u)$ et $\mathfrak{E}^-:=\star^{\dot{\phi}_{t\to- \infty}}(u)$ tels que 
\begin{equation}
\mathcal{O}(\mathfrak{E}^+):=\mathcal{O}^{\dot{\phi}_{t\to + \infty}}(u) \ \ \ \textnormal{et} \ \ \ \mathcal{O}(\mathfrak{E}^-):=\mathcal{O}^{\dot{\phi}_{t\to - \infty}}(v).
\end{equation}

\section{Preuve complémentaire}
\begin{lemma}
Soient $S_1$ et $S_2$ deux sous-variétés de la $n$-variété $M$. Si $S_1\pitchfork S_2$, alors $\forall \ z\in S_1\cap S_2$
\begin{equation}
\mathcal{T}_z(S_1\cap S_2)=\mathcal{T}_zS_1\cap \mathcal{T}_zS_2.
\end{equation}
\end{lemma}
\begin{prL}
Comme $S_1\cap S_2\subseteq S_1$ et $S_1\cap S_2\subseteq S_2$, on a pour tout $z\in S_1\cap S_2$
\begin{equation}
\mathcal{T}_z(S_1\cap S_2)\subseteq \mathcal{T}_zS_1 \ \textnormal{et} \ \mathcal{T}_z(S_1\cap S_2)\subseteq \mathcal{T}_zS_2 \Rightarrow \mathcal{T}_z(S_1\cap S_2)\subseteq \mathcal{T}_zS_1 \cap \mathcal{T}_zS_2.
\end{equation}
Nonobstant cela, la condition de transversalité sur l'intersection  implique que 
\begin{equation}
\textnormal{codim}(S_1\cap S_2)=\textnormal{codim}(S_1)+\textnormal{codim}(S_2),
\end{equation}
qui se traduit par
\begin{equation}
\textnormal{dim}(S_1\cap S_2)=\textnormal{dim}(S_1)+\textnormal{dim}(S_2)-n \Rightarrow \textnormal{dim}(\mathcal{T}_z(S_1\cap S_2))=\textnormal{dim}(\mathcal{T}_zS_1)+\textnormal{dim}(\mathcal{T}_zS_2)-n.
\end{equation}
La transversalité implique aussi que 
\begin{equation}
\textnormal{dim}(\mathcal{T}_zM)=\textnormal{dim}(\mathcal{T}_zS_1 + \mathcal{T}_zS_2)=n.
\end{equation}
Nous réécrivons alors 
\begin{equation}
\textnormal{dim}(\mathcal{T}_z(S_1\cap S_2))=\textnormal{dim}(\mathcal{T}_zS_1)+\textnormal{dim}(\mathcal{T}_zS_2)-\textnormal{dim}(\mathcal{T}_zS_1 + \mathcal{T}_zS_2).
\end{equation}
du théorème fondamental de l'algèbre linéaire, 
\begin{equation}
\textnormal{dim}(\mathcal{T}_zS_1\cap \mathcal{T}_zS_2)=\textnormal{dim}(\mathcal{T}_zS_1)+\textnormal{dim}(\mathcal{T}_zS_2)-\textnormal{dim}(\mathcal{T}_zS_1 + \mathcal{T}_zS_2).
\end{equation}
Ainsi $\textnormal{dim}(\mathcal{T}_zS_1\cap \mathcal{T}_zS_2)=\textnormal{dim}(\mathcal{T}_z(S_1\cap S_2))$ et $\mathcal{T}_z(S_1\cap S_2)\subseteq \mathcal{T}_zS_1 \cap \mathcal{T}_zS_2$ implique le résultat désiré.
\end{prL}







\bibliography{ete2018ref}
	\bibliographystyle{acm}
\end{document}